\documentclass[12pt]{article}
\usepackage{amsthm}
\usepackage{amsfonts}
\usepackage{amsmath}
\usepackage{amssymb}
\usepackage{hyperref}
\usepackage{epsfig}
\newcommand{\crn}{\ensuremath{\mathrm{cr}}}
\newtheorem{theorem}{Theorem}
\newtheorem{observation}[theorem]{Observation}

\newtheorem{megaclaim}{Claim}
\newtheorem{problem}{Problem}
\newcommand{\claim}[2]{\begin{megaclaim}\label{#1}  #2 \end{megaclaim}}
\newcommand{\cref}[1]{Claim~\ref{#1}}

\def\myincludegraphics#1{\begin{center}\includegraphics{#1}\end{center}}

\title{Graphs with two crossings are $5$-choosable\thanks{Supported by the CZ-SL bilateral project MEB 091037 and BI-CZ/10-11-004 and grants
GACR~201/09/0197 and 1M0021620808. 
}}
\author{Zden\v{e}k Dvo\v{r}\'ak\thanks{Charles University, Prague, Czech Republic. E-mail: {\tt rakdver@kam.mff.cuni.cz}.}\and
Bernad Lidick\'y\thanks{Charles University, Prague, Czech Republic. E-mail: {\tt bernard@kam.mff.cuni.cz}. Supported by GAUK 60310.}\and
Riste \v{S}krekovski\thanks{University of Ljubljana, Ljubljana, Slovenia. E-mail: {\tt skrekovski@gmail.com}.}}
\date{}

\begin{document}
\maketitle

\begin{abstract}
A graph $G$ is $k$-choosable if $G$ can be properly colored
whenever every vertex has a list of at least $k$ available colors.
Thomassen's theorem
states that every planar graph is $5$-choosable.
We extend the result by showing that every graph with 
at most two crossings is $5$-choosable.
\end{abstract}

\section{Introduction}

All graphs considered in this paper are simple, i.e., without loops or parallel edges.
We denote the set of vertices of a graph $G$ by $V(G)$ and the set of edges by $E(G)$.
The \emph{crossing number} of $G$, denoted by $\crn(G)$, is the minimum possible number of crossings 
in a drawing of $G$ in the plane.  


Let $G$ be a graph and $C$ a set of colors.  A {\em list assignment} is a function $L:V(G) \to 2^C$.
We say that $G$ is {\em $L$-colorable} if there exists a coloring $\varphi: V(G) \to C$
such that $\varphi(v) \in L(v)$ for every vertex $v\in V(G)$ and adjacent vertices are assigned
different colors. We say that a graph $G$ is {\em $k$-choosable} if $G$ is $L$-colorable
whenever $L$ assigns at least $k$ colors to each vertex.

The concept of list colorings and choosability was introduced by Vizing~\cite{vizing1976} and independently by 
Erd\H{o}s et al.~\cite{erdosrubintaylor1979}.  Clearly, if a graph is $k$-choosable, then it is also $k$-colorable.
However, the opposite implication does not hold.  For instance, there exist planar graphs that are not $4$-choosable,
see Voigt~\cite{voigt1993}.  On the other hand, Thomassen~\cite{thomassen1994} gave a strikingly beautiful proof that
every planar graph is $5$-choosable.

A general bound for choosability of graphs on surfaces,
an analogue of Dirac's Map-Color Theorem, is known due to B\"ohme et al.~\cite{BoMoSt99}.  However, only very
few graphs have choosability close to this bound.  DeVoss et al.~\cite{DeKaMo08} obtained more general result claiming
that locally planar graphs are 5-choosable, which extended the result of Thomassen to non-planar graphs.
Furthermore, Kawarabayashi and Mohar~\cite{kawmoh} proved that there are only finitely 
many minimal graphs that are not $5$-choosable on any fixed surface.

Note that the genus of any graph is bounded by its crossing number, thus it might be possible to obtain more refined results
by considering graphs with bounded crossing number.  For example, for ordinary chromatic number it is easy to see that
all graphs with crossing number at most two are $5$-colorable.  Oporowski and Zhao~\cite{opzha} improved this result
by showing that $K_6$ is the only $6$-critical graph with crossing number three, and conjectured that it is the only
$6$-critical graph with crossing number at most $5$.  This was eventually settled by Erman et al.~\cite{EHLP-11},
who proved that $K_6$ is the only $6$-critical graph with crossing number at most $4$ and gave an example of a $6$-critical
graph $K_6$-free graph with crossing number $5$.

Similarly, the result of Thomassen can be easily used to derive that every graph with
at most one crossing is 5-choosable.  Erman et al.~\cite{EHLP-11} posed a question if this is also
true for graphs with two crossings.  The main result of this paper is the positive answer to this question:

\begin{theorem}\label{thm-main}
Every graph $G$ with $\crn(G) \leq 2$ is 5-choosable.
\end{theorem}

Inspired by the result of Erman et al.~\cite{EHLP-11}, we pose the following open problem:
\begin{problem}
Is it true that every $K_6$-free graph $G$ with $\crn(G) \leq 4$ is $5$-choosable?
\end{problem}

Note that the Four Color Theorem is used heavily in the coloring case, thus it is not obvious that
the results should generalize to the choosability case.

\section{$5$-choosability of graphs with two crossings} 

In order to show $5$-choosability of planar graphs, Thomassen~\cite{thomassen1994} proved
the following more general statement.

\begin{theorem}\label{thm-thom}
Let $G$ be a plane graph, $F$ a face of $G$ and $xy$ an edge incident with $F$.  Then $G$ is $L$-colorable
for any list assignment $L$ such that
\begin{itemize}
\item $|L(v)|\ge 5$ for $v\in V(G)\setminus V(F)$,
\item $|L(v)|\ge 3$ for $v\in V(F)\setminus\{x,y\}$,
\item $|L(x)|,|L(y)|\ge 1$, and
\item if $|L(x)|=|L(y)|=1$, then $L(x)\neq L(y)$.
\end{itemize}
\end{theorem}

We frequently use the following observation:
\begin{observation}\label{cor-drtic}
Let $G$ be a connected graph with a list assignment $L$.  Suppose that the subgraph $G[S]$ of $G$ induced by a set
$S \subset V(G)$ has an $L$-coloring $\psi$.
Let $G'=G-S$ and for every $v \in V(G')$ let $L'(v) = L(v) \setminus \{\psi(u): u \in N(v) \cap S\}$.
If $G'$ and $L'$ satisfy the assumptions of Theorem~\ref{thm-thom} then
$\psi$ can be extended to a coloring of $G$.
\end{observation}

Instead of proving Theorem~\ref{thm-main} directly, we prove a slightly stronger theorem.
\begin{theorem}\label{thm-prove}
Let $G$ be a graph and $L$ a list assignment such that either
\begin{itemize}
\item $\crn(G) \leq 2$ and $|L(v)|\ge 5$ for every $v\in V(G)$, or
\item $\crn(G) \leq 1$, $G$ contains a triangle $T$, $L(v) = 1$ for all $v\in V(T)$,
$L(u)\neq L(v)$ if $u$ and $v$ are two distinct vertices of $T$ and
$|L(v)| \geq 5$ for all $v \in V(G)\setminus V(T)$.
\end{itemize}
Then $G$ is $L$-choosable.
\end{theorem}

\begin{proof}[Proof of Theorem~\ref{thm-prove}]

Let $G$ and $L$ be a counterexample with the smallest crossing number, subject to that the smallest number of vertices,
and subject to that the largest number of edges. 
Observe that $G$ is connected, as otherwise we can color each connected
component of $G$ separately.  Furthermore, the minimum degree of $v \in V(G)\setminus V(T)$ is $5$, 
as if $v$ had degree at most $4$, then an $L$-coloring of $G-v$ 
(which exists by the minimality of $G$) can be extended to an $L$-coloring of $G$.
Moreover, we assume that $|L(v)| = 5$ for all $v \in V(G)\setminus V(T)$ as removing
colors from lists does not turn $G$ into an $L$-colorable graph.

\claim{claim:basic}{$\crn(G) \geq 1$.}
\begin{proof}
Suppose for contradiction that $G$ is planar. By Theorem~\ref{thm-thom}, it contains
a precolored triangle $T=t_1t_2t_3$. Let $S=\{t_1\}$ and let $\psi$ be a coloring of $S$ such that $\psi(t_1) \in L(t_1)$. 
By applying Observation~\ref{cor-drtic} we get an $L$-coloring of $G$.
\end{proof}

We call crossings and $T$ \emph{dangerous configurations}.
Since all graphs with crossing number one are $5$-choosable by~\cite{EHLP-11},
we can assume that $G$ has two dangerous configurations: either $\crn(G)=2$ or
$\crn(G)=1$ and $G$ contains $T$.  Let us fix a drawing of $G$ with the minimum number of crossings.

\claim{cl:noct}{If $T$ exists, then no edge of $T$ is crossed.}
\begin{proof}
Let $T=uvw$, and assume for contradiction that the edge $uv$ is crossed by an
edge $xy$.  Let $G_1$ be the subgraph of $G$ induced by the vertices
drawn in the closed disk bounded by $T$, and let $G_2=G-(V(G_1)\setminus V(T))$.
By symmetry, assume that $x\neq w$ and $x\in V(G_1)$.  There exists an $L$-coloring
of $G_2$, and by Observation~\ref{cor-drtic} applied with $S=V(G_2)\setminus \{v,w\}$,
this coloring extends to an $L$-coloring of $G$.  This is a contradiction.
\end{proof}

Let $C$ be a cycle in $G$.  Let $G_1$ be the subgraph of $G$ consisting of vertices and edges
drawn in the \emph{closed} disk bounded by $C$ and $G_2$ the subgraph of $G$ consisting of the vertices
and edges drawn outside of the \emph{open} disk bounded by $C$.  Note that $G_1\cap G_2=C$.
If no edge of $C$ is crossed and $V(G_1)\neq V(C)\neq V(G_2)$, then we say that $C$ is a \emph{separating cycle}.
We call $G_1$ and $G_2$ the \emph{$C$-components}.

\claim{claim:triangle}{$G$ does not contain a separating triangle.}
\begin{proof}
Suppose for a contradiction that there is a separating triangle $C = x_1x_2x_3$,
and let $G_1$ and $G_2$ be the $C$-components.
If both dangerous configurations are in $G_1$, we first color $G_1$ by induction
and then extend coloring of $C$ to $G_2$ by \cref{claim:basic}.
Otherwise, without loss of generality, we assume that if $\crn(G)=1$, then $T\subset G_1$.  
We first color $G_1$ and then extend the coloring of $C$ to $G_2$, where $C$ plays the
role of the precolored triangle in $G_2$.
\end{proof}

Similarly (by adding vertices to extend the cut to a triangle if necessary) one can prove the following.
\claim{claim:conn}{$G$ is $2$-connected and if $\{u,v\}$ is a cut in $G$, then $uv$ is not a non-crossed edge.}

Furthermore, we can restrict separating $4$-cycles.
\claim{claim:sep}{$G$ does not contain a separating cycle $C$ of length four with both dangerous configurations
draw on the same side of $C$.}
\begin{proof}
Suppose for a contradiction that there is a separating cycle $C = x_1x_2x_3x_4$
with $C$-components $G_1$ and $G_2$, where both dangerous configurations are in $G_1$.
By \cref{claim:triangle}, we can assume that $C$ is an induced cycle in $G_2$.
By the minimality of $G$, there exists an $L$-coloring $\psi$ of
$G_1$.  Observation~\ref{cor-drtic} used with $S=V(G_1)\setminus \{x_1,x_2\}$
implies that $\psi$ can be extended to an $L$-coloring of $G$, which is a contradiction.
\end{proof}

\claim{claim:twicecross}{No edge is crossed twice.}
\begin{proof}
Suppose for a contradiction that an edge $e=uv$ is crossed by edges $e_1$ and $e_2$.
We distinguish two cases depending on the number of vertices incident with $e_1$ and $e_2$.

Suppose first that there exists a vertex $w$ incident to both $e_1$ and $e_2$. 
Let $S=\{u,w\}$ and let $\psi$ be an arbitrary $L$-coloring of $S$.
Observation~\ref{cor-drtic} implies that $\psi$ can be extended to an $L$-coloring of $G$, which is a contradiction.

Therefore, no vertex is incident with both $e_1$ and $e_2$.  Let $e_1=w_1z_1$ and $e_2=w_2z_2$.
As $G$ has the largest possible number of edges, we can assume that $C=uw_1w_2vz_2z_1$ is a cycle of
length $6$ with both crossings drawn inside.  By \cref{claim:sep}, no vertex is drawn inside $C$.
Let $S=\{u,w_1,w_2\}$ and let $B$ be the set of common neighbors of the vertices of $S$.

Since every vertex of $C$ has degree at least $5$, Claims~\ref{claim:triangle}, \ref{claim:conn} and \ref{claim:sep}
imply that $B \cap V(C) = \emptyset$. Suppose there exists $x \in B$. \cref{claim:sep} implies that
triangles $uw_1x$ and $xw_1w_2$ bound faces. Hence $w_1$ has degree four, contradicting that the minimum degree of $G$ is $5$. So, $B = \emptyset$.
We conclude that we can apply Observation~\ref{cor-drtic} for an arbitrary $L$-coloring of $G[S]$
and obtain an $L$-coloring of $G$. 
\end{proof}

It turns out that we can restrict our attention only to the case that $T$ exists.
\claim{claim:cr1}{$\crn(G)=1$.}
\begin{proof}
For contradiction, assume that $\crn(G)=2$, and
let edges $e=xx'$ and $f=yy'$ cross each other. Let $X=\{x,x',y,y'\}$.
Let $G'$ be the graph obtained
from $G-\{e,f\}$ by adding a new vertex $v$ adjacent to all vertices of $X$.
Note that $xy$ is an edge as $G$ has the largest number of edges. 
Let $L'$ be the list assignment such that $L'(x)\subseteq L(x)$ and $L'(y)\subseteq L(y)$
are distinct lists of size one, $L'(v)=\{c\}$ for a new color $c$ that does not appear in any
of the lists, $L'(x')=(L(x')\setminus L'(x))\cup \{c\}$, $L'(y')=(L(y')\setminus L'(y))\cup \{c\}$
and $L'(z)=L(z)$ for every $z\in V(G)\setminus X$.
Since $\crn(G')<\crn(G)$, the graph $G'$ (with
$xyv$ playing the role of the precolored triangle) has an $L'$-coloring $\varphi$.  Note that
$\varphi(x)\neq\varphi(x')\neq c$ and $\varphi(y)\neq\varphi(y')\neq c$, hence $\varphi$ is also
an $L$-coloring of $G$.
\end{proof}

Let $X =\{v_1,v_2,v_3,v_4\}$, where $e=v_1v_3$ and $f=v_2v_4$ cross each other. 
Since $G$ has the largest possible number of edges, the following claim holds.

\claim{cl-fk4}{$G[X]$ is a complete graph.}  

So, $v_1v_2v_3v_4$ is a cycle of length four enumerated in the clockwise order.
Let the precolored triangle $T$ have vertices $t_1$, $t_2$ and $t_3$ in the clockwise order.

\claim{claim:noboth}{Let $h=uv$ be equal to $e$ or $f$.  If a vertex $w$ is adjacent
to both $u$ and $v$, then $w$ belongs to $X$.}
\begin{proof}
Suppose for a contradiction that $w$ does not belong to $X$ and let $\{x,y\}=X\setminus \{u,v\}$.
By symmetry between the cycles $wuxv$ and $wuyv$, we can assume that both dangerous configurations appear inside the
closed disk bounded by $wuxv$, and by \cref{claim:sep} the cycle $wuxv$ is not separating.  We conclude
that $x$ has degree at most four, which is a contradiction.
\end{proof}

\claim{claim:sepcross}{$G$ does not contain a separating cycle $C$ of length four
such that $|V(C)\cap V(X)|\ge 2$.}
\begin{proof}
Suppose for a contradiction that $C$ is such a cycle.
Note that no edge of $C$ is crossed as $C$ is separating.
Let $G_1$ and $G_2$ be the $C$-components, where $X\subseteq V(G_2)$, and let $u$ and $v$ be two vertices in $V(C)\cap X$.
\cref{claim:sep} implies that both $G_1$ and $G_2$ contain a dangerous configuration.

By \cref{cl-fk4}, $u$ and $v$ are adjacent in $G_2$.  By \cref{claim:noboth}, $uv$ is not a crossed edge.
By \cref{claim:triangle}, we conclude that $uv$ is an edge of $C$ and that $C$ is an induced cycle in $G_2$.
By the minimality of $G$, there exists an $L$-coloring $\psi$ of $G_1$.  
Observation~\ref{cor-drtic} used with $S=(V(G_1)\setminus V(C))\cup \{u,v\}$
implies that $\psi$ can be extended to an $L$-coloring of $G$, which is a contradiction.
\end{proof}

\claim{claim:indep}{$V(T)$ and $X$ are disjoint.}
\begin{proof}
Let $uv$ be a non-crossed edge such that $u \in V(T) \cap X$
and $v \in X$. Let $S = \{u,v\}$
and $\psi$ be an arbitrary $L$-coloring of $G[S\cup V(T)]$. See Figure~\ref{fig-onedge}(a).
Observe that $G - S$ is planar and all neighbors of $S$ are incident with one face. 
Hence, we can apply Observation~\ref{cor-drtic} and obtain an $L$-coloring of $G$, which is a contradiction.
\end{proof}


\begin{figure}
\myincludegraphics{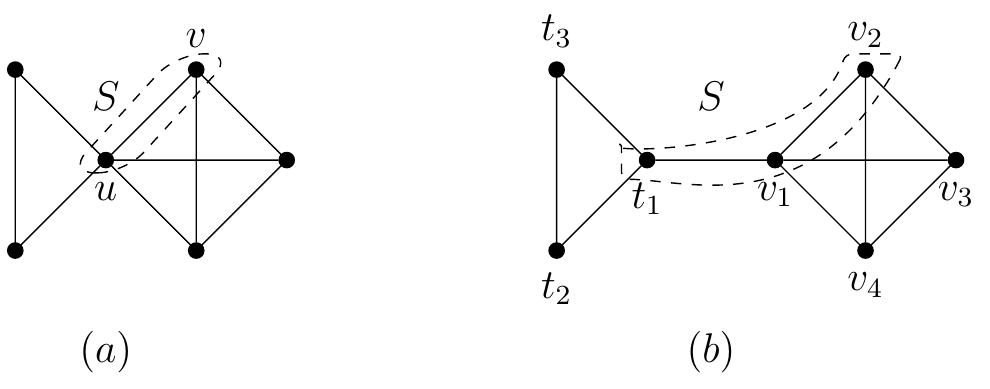}
\caption{The configurations of Claims \ref{claim:indep} and \ref{claim-dist1}.\label{fig-onedge}}
\end{figure}

\claim{claim-dist1}{There is no edge $tv$ such that $t \in V(T)$ and $v \in X$.}
\begin{proof}
Suppose without loss of generality that $t_1v_1$ is an edge.
Let $S=\{t_1,v_1,v_2\}$, see Figure~\ref{fig-onedge}(b).
By symmetry between $v_2$ and $v_4$ and by \cref{claim:noboth}, we can assume that $t_1v_4\not\in E(G)$.

Let us construct an $L$-coloring $\psi$ of $S$.
Let $\psi(t_1)$ be the unique color in $L(t_1)$ and choose the color of $v_1$ so that
$\psi(v_1) \in L(v_1) \setminus( L(t_1) \cup L(t_2) \cup L(t_3))$.
Now we need to choose $\psi(v_2)$ such that all vertices of $G \setminus S$
except for $t_2$ and $t_3$ have at least three colors left in their lists. 
Hence we only need to deal with vertices in $V(G)\setminus V(T)$ adjacent to all 
three of $t_1$, $v_1$ and $v_2$. We call such vertex $y$ a \emph{troublemaker}
if $\psi(t_1),\psi(v_1) \in L(y)$.
If there is no troublemaker, then we choose $\psi(v_2) \in L(v_2) \setminus (\{\psi(v_1), \psi(t_1)\} \cup L(t_2) \cup L(t_3))$
and use Observation~\ref{cor-drtic} to obtain an $L$-coloring of $G$.

Since $t_1v_4\not\in E(G)$, $v_4$ is not a troublemaker.  Furthermore, if $t_1v_3 \in E(G)$, then $v_4$ or $v_2$ would have degree
at most four by \cref{claim:sep}, which is a contradiction.  Consequently, $v_3$ is not a troublemaker.
Suppose that there is a troublemaker $y \in V(G)\setminus (X \cup V(T))$.  By \cref{claim:triangle},
$t_1v_1y$ and $v_1v_2y$ are faces (hence, there is no other troublemaker) and $t_1$ is not adjacent to $v_2$.
\cref{claim:sepcross} implies that $v_2$ is adjacent to neither $t_2$ nor $t_3$.

We choose $\psi(v_2)$ arbitrarily from $L(v_2) \setminus (L(y) \setminus \{\psi(t_1)\})$.
Note that there is at least one choice for $\psi(v_2)$, since $\psi(t_1)\in L(y)$.
Because $\psi(v_1) \in L(y)$, the resulting coloring of $G[S]$ is proper.
Furthermore, $|L(y)\setminus \{\psi(t_1),\psi(v_1),\psi(v_2)\}|=3$, hence
Observation~\ref{cor-drtic} applies.
\end{proof}

Let $P = p_1\ldots p_k$ be a path such that $p_1 \in V(T)$
and $p_{k-1},p_k \in X$ and no edge of $P$ is crossed, see Figure~\ref{fig-p}.
Let the {\em score} of this path be $2k-b$, where 
$$b = \begin{cases}  1  & \text{if\ } p_{k-2}p_k\in E(G) \\  0  & \text{otherwise.} \end{cases}$$ 

Let $P\subseteq G$ be such a path with the smallest possible score.
By \cref{claim-dist1}, we have $k\ge 4$.  Let $\overline{P}=p_1p_2\ldots p_{k-1}$, see Figure~\ref{fig-p}.
The path $\overline{P}$ is induced, as otherwise $P$ would contain a shorter subpath
with a smaller score.  Assume without loss of generality that $p_1 = t_1$.

\begin{figure}
\myincludegraphics{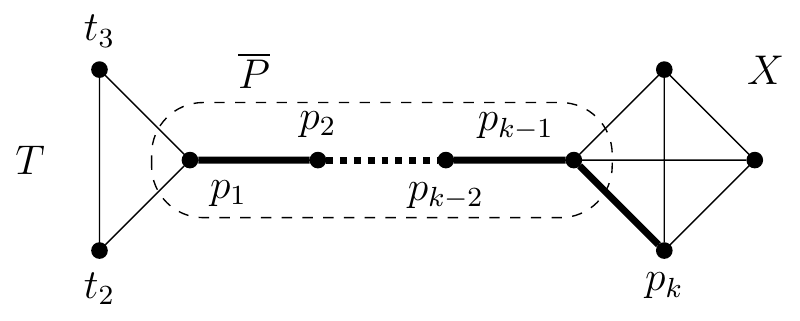}
\caption{Paths $P$ and $\overline{P}$ connecting dangerous configurations.\label{fig-p}}
\end{figure}

\claim{claim:adj}{If $p_i,p_j\in V(\overline{P})$ are neighbors of a vertex $v \in V(G)\setminus V(P)$, then $|i-j|\le 2$.}
\begin{proof}
Suppose that $i\ge j+3$.  Then the path $p_1\ldots p_jvp_i\ldots p_k$ is shorter than $P$,
contradicting the choice of $P$.  
\end{proof}

\claim{claim:three}{Every vertex $v \in V(G)\setminus V(P)$ has at most three neighbors in $P$.}
\begin{proof}
Suppose for contradiction that $v$ has at least four neighbors in $P$. \cref{claim:adj} implies
that $v$ has four neighbors and they are $p_{k-3},p_{k-2},p_{k-1}$ and $p_{k}$.
Since the path $p_1\ldots p_{k-3}vp_kp_{k-1}$ does not have smaller score than $P$,
it follows that $p_k$ is adjacent to $p_{k-2}$.
But that gives a contradiction with \cref{claim:triangle}, since either $p_{k-1}p_kv$ or $p_{k-2}p_{k-1}p_k$ is separating.
\end{proof}

Observe that $t_2$ and $t_3$ have at most two neighbors in $P$.  Furthermore, consider a vertex $u\in X$; if $u$ has three neighbors
in $P$, then by the choice of $P$, these neighbors are $p_k$, $p_{k-1}$ and $p_{k-2}$.  By \cref{claim:noboth}, $p_{k-2}$ is not
adjacent to $p_k$.  This contradicts the choice of $P$, since the path $p_1\ldots p_{k-1}u$ has smaller score.
Therefore, we have the following.
\claim{claim:notinx}{If $v\in V(G)\setminus V(P)$ has three neighbors in $P$, then $v\not\in V(T)\cup X$.}

For every vertex $v\in V(G)\setminus V(P)$ with three neighbors in $P$, let $g_P(v)=p_i$, where $p_i\in V(P)$ is the neighbor of $v$
with the largest $i$.  We define $g_P(v)=v$ if $v$ has at most two neighbors in $P$.  We write $g(v)$ instead of $g_P(v)$ for brevity when
the path $P$ is clear from the context.

\claim{claim:uniq}{If $u$ and $v$ are distinct vertices of $V(G)\setminus V(P)$, then $g(u) \neq g(v)$.}
\begin{proof}
Let $p_g=g(u)=g(v)$ for two distinct vertices $u$ and $v$.
If $g\neq k$, then both $u$ and $v$ are adjacent to $p_{g-2}$, $p_{g-1}$ and $p_g$ by \cref{claim:adj}.  However, that contradicts
Claim \ref{claim:triangle} or \ref{claim:sep}.  Hence, we have $g=k$.

By the choice of $P$, all neighbors of $u$ and $v$ in $P$ are contained
in $\{p_{k-3}, p_{k-2}, p_{k-1}, p_k\}$.  By \cref{claim:triangle}, $u$ and $v$ cannot both be adjacent to $p_{k-1}$,
thus assume that say $u$ is adjacent to $p_{k-3}$, $p_{k-2}$ and $p_k$.  By Claims \ref{claim:sep} and \ref{claim:sepcross},
the cycle $p_{k-2}p_{k-1}p_ku$ is not separating, hence by \cref{claim:triangle}, $v$ is not adjacent to $p_{k-1}$.
But then $v$ is adjacent to $p_{k-2}$ and the cycle $p_{k-2}p_{k-1}p_kv$ is separating, which is a contradiction.
\end{proof}


Let $S=V(P)$.  We now attempt to construct an $L$-coloring $\psi$ of $G[S]$ so that the assumptions of Observation~\ref{cor-drtic}
are satisfied.  

We will assign the colors to all vertices of $P$ in order.  We start with the unique choice $\psi(p_1) \in L(p_1)$
and color $p_2$ by a color $\psi(p_2)\in L(p_2)\setminus (L(t_1)\cup L(t_2)\cup L(t_3))$.  Note that no other
vertex of $P$ has a neighbor in $T$.
Suppose that we have already colored the vertices $p_1$, \ldots, $p_{j-1}$ and let $R_j=V(G)\setminus \{t_2,t_3,p_1,\ldots, p_{j-1}\}$.  For a vertex $v\in R_j$,
let $B_j(v)=L(v)\setminus \{\psi(p_i) : 1\le i\le j-1, vp_i\in E(G)\}$.  
We choose the color $\psi(p_j)\in B_j(p_j)$ in such a way that $|B_{j+1}(v)|\ge 3$ for any $v\in R_{j+1}$.
This coloring of $P$ ensures that all vertices of $G-S$ other than $t_2$ and $t_3$ have at least three available colors; and since
$t_2$ and $t_3$ are adjacent
we can apply Observation~\ref{cor-drtic} and obtain an $L$-coloring of $G$, which is a contradiction.

Let us now describe how $\psi(p_j)$ is chosen.  Let $y$ be a vertex such that $g(y)=p_j$ if such a vertex exists.
Regardless of the choice of $\psi(p_j)$, for any vertex $v\in R_{j+1}$ other than $y$ we have $|B_{j+1}(v)|\ge 3$, since $v$ has at most two neighbors
in $\{p_1,\ldots, p_j\}$ or it is not adjacent to $p_j$.
The same holds for $y$ if $|B_j(y)|\ge 4$, thus assume that $|B_j(y)|=3$.  
If $B_j(p_j)\not\subseteq B_j(y)$, then we can choose
$\psi(p_j)\in B_j(p_j)\setminus B_j(y)$.  Therefore, since $|B_j(p_j)|\ge 3$, we can assume that $B_j(p_j)=B_j(y)$.  Since $|B_j(p_j)|=3$, $p_j$ has two neighbors $p_i, p_l\in V(P)$
with $i<l<j$.
Since the path $\overline{P}$ is induced, this is only possible if $j=k$ and $p_k$ is adjacent to $p_{k-2}$.
By \cref{claim:triangle}, $y$ is not adjacent to $p_{k-1}$. 
Since $y$ has exactly three neighbors in $P$, $y$ is adjacent to $p_{k-3}$, $p_{k-2}$ and $p_k$.

Consider now the path $P' = p_1\ldots p_{k-2}p_{k}p_{k-1}$ instead of $P$.  Note that $P'$ has the same score as $P$,
thus we conclude that there also exists a vertex $y'$ adjacent to $p_{k-1}$, $p_{k-2}$ and $p_{k-3}$.
By \cref{claim:triangle}, $yp_{k-3}p_{k-2}$, $yp_{k-2}p_k$, $y'p_{k-3}p_{k-2}$ and $y'p_{k-2}p_{k-1}$ are
faces, see Figure~\ref{fig-K}.

\begin{figure}
\myincludegraphics{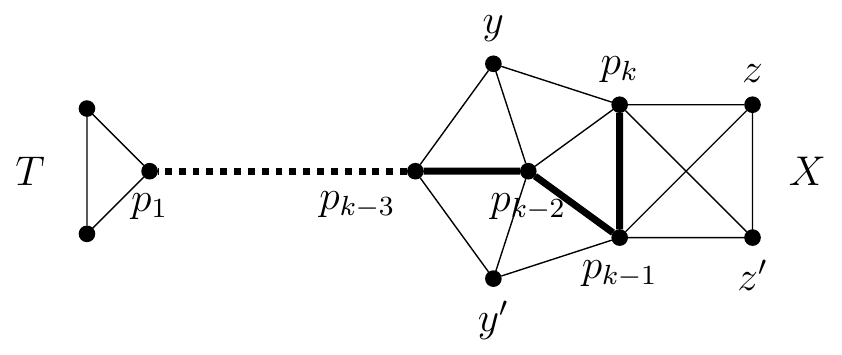}
\caption{Positions of $T$, $X$, $y$ and $y'$.\label{fig-K}}
\end{figure}

Let us fix a coloring $\psi$ of $P-\{p_{k-2},p_{k-1},p_k\}$ such that $|B_{k-2}(v)|\ge 3$ for every $v\in R_{k-2}$.  
By the preceding arguments, we can assume the following.
\claim{cl-nex}{There is no $L$-coloring of $G[V(T)\cup V(P-p_k)]$ such that $|B_k(v)|\ge 3$ for every $v\in R_k$
and $|B_k(y)|\ge 4$ or $B_k(p_k)\neq B_k(y)$.}  

The only neighbors of $p_{k-2}$ in $V(G)\setminus V(P)$ are $y$ and $y'$, hence there is no vertex $v\in V(G) \setminus V(P)$ with $g(v)=p_{k-2}$.
Since $\psi$ cannot be extended to a coloring contradicting \cref{cl-nex}, we will show that:

\claim{claim:bsame}{$B_{k-2}(p_{k-2}) = B_{k-2}(y) = B_{k-2}(y')$ and 
$L(p_{k-1}) = L(p_{k}) = B_{k-2}(p_{k-2}) \cup \{c\}$ for some color $c$.}
\begin{proof}
Since $p_{k-3}$ is the only neighbor of $y$, $y'$ and $p_{k-2}$ in $P-\{p_{k-2},p_{k-1},p_k\}$, we have
$|B_{k-2}(y)|,|B_{k-2}(y')|,|B_{k-2}(p_{k-2})| \geq 4$.  If $|B_{k-2}(y)|=5$ or $B_{k-2}(p_{k-2})\not\subseteq B_{k-2}(y)$,
we can choose $\psi(p_{k-2})\in B_{k-2}(p_{k-2})$ so that $|B_{k-1}(y)|\ge 4$, and further extend $\psi$ to a coloring
of $P-p_k$ contradicting \cref{cl-nex}.  It follows that $|B_{k-2}(y)|=4$ and $B_{k-2}(p_{k-2}) = B_{k-2}(y)$.
Symmetrically, we obtain $B_{k-2}(p_{k-2}) = B_{k-2}(y')$ by considering $P'$.

Consider colors $c_1 \in L(p_{k-1}) \setminus B_{k-2}(p_{k-2})$ and $c_2 \in L(p_{k}) \setminus B_{k-2}(p_{k-2})$.
If $c_1\neq c_2$, then choose $\psi(p_{k-2})\in B_{k-2}(p_{k-2})$ arbitrarily and set $\psi(p_{k-1}) = c_2$.
Since $c_2\in B_k(p_k)\setminus B_k(y)$, this coloring contradicts \cref{cl-nex}.  Therefore, $c_1=c_2$,
which implies that $L(p_{k-1}) = L(p_{k}) = B_{k-2}(p_{k-2}) \cup \{c_1\}$.
\end{proof}

\begin{figure}
\myincludegraphics{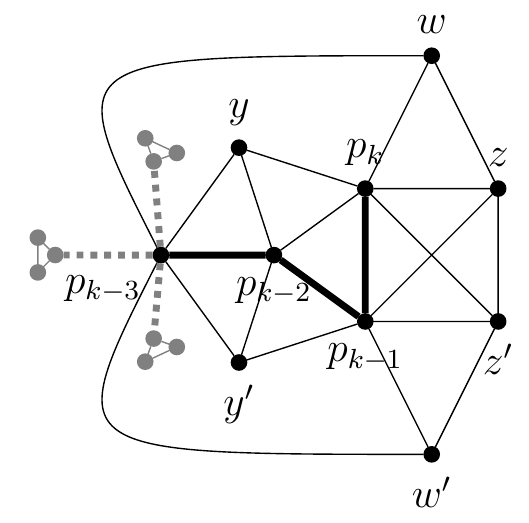}
\caption{Vertices $w$ and $w'$ and position of $T$.\label{fig-w}}
\end{figure}

Let us now choose $\psi(p_{k-2})\in B_{k-2}(p_{k-2})$ arbitrarily and set $\psi(p_{k-1})=c$.
Note that $B_k(y)=B_k(y')=B_k(p_k)$.

Let $X\setminus \{p_{k-1},p_k\}=\{z,z'\}$, where $z$ is joined to $p_k$ by a non-crossed edge.
By planarity, $z'$ is not adjacent to $p_{k-2}$, thus $B_k(z') \ge 4$.
Let $Q=p_1\ldots p_{k-1}z'$.  If no vertex $w'\not\in\{y,y',p_k\}$ satisfies $g_Q(w')=z'$,
then we can choose $\psi(z')\in B_k(z')\setminus B_k(y')$ and apply Observation~\ref{cor-drtic}
with $S=V(Q)$, obtaining an $L$-coloring of $G$.  Therefore, we may assume that there
exists such a vertex $w'$.  Since $w'$ has at least three neighbors in $Q$ and it is not adjacent to $p_{k-2}$ by planarity,
the choice of $P$ implies that $w'$ is adjacent to $p_{k-1}$ and $p_{k-3}$.  Symmetrically, by considering
the path $Q'=p_1\ldots p_{k-2}p_kz$, we conclude that there exists a vertex $w\not\in\{y,y',p_{k-1}\}$
adjacent to $z$, $p_k$ and $p_{k-3}$.  However, by planarity either $p_{k-3}p_{k-2}p_{k-1}w'$ or $p_{k-3}p_{k-2}p_kw$
contradicts \cref{claim:sep}, see Figure~\ref{fig-w} showing the possible positions of $T$ with respect to these cycles.
\end{proof}

\bibliographystyle{siam}
\bibliography{2cross}

\begin{thebibliography}{1}

\bibitem{BoMoSt99}
{\sc T.~B\"ohme, B.~Mohar, and M.~Stiebitz}, {\em Dirac's map-color theorem for
  choosability}, Journal of Graph Theory, 32 (1999), pp.~311--326.

\bibitem{DeKaMo08}
{\sc M.~DeVos, K.~Kawarabayashi, and B.~Mohar}, {\em Locally planar graphs are
  5-choosable}, Journal of Combinatorial Theory, Series B, 98 (2008),
  pp.~1215--1232.

\bibitem{erdosrubintaylor1979}
{\sc P.~Erd\H{o}s, A.~L. Rubin, and H.~Taylor}, {\em Choosability in graphs},
  Congressus Numerantium, 26 (1980), pp.~125--157.

\bibitem{EHLP-11}
{\sc R.~Erman, F.~Havet, B.~Lidick\'{y}, and O.~Pangr\'{a}c}, {\em 5-colouring
  graphs with 4 crossings}, SIAM Journal on Discrete Mathematics,  (2010).
\newblock to appear.

\bibitem{kawmoh}
{\sc K.~Kawarabayashi and B.~Mohar}, {\em List-color-critical graphs on a fixed
  surface}, in Proceedings of the Twentieth Annual ACM-SIAM Symposium on
  Discrete Algorithms (SODA'09), SIAM, 2009, pp.~1156--1165.

\bibitem{opzha}
{\sc B.~Oporowski and D.~Zhao}, {\em Coloring graphs with crossings}, Discrete
  Mathematics, 309 (2009), pp.~2948--2951.

\bibitem{thomassen1994}
{\sc C.~Thomassen}, {\em Every planar graph is 5-choosable}, Journal of
  Combinatorial Theory, Series B, 62 (1994), pp.~180--181.

\bibitem{vizing1976}
{\sc V.~G. Vizing}, {\em Vertex colorings with given colors (in russian)},
  Metody Diskretnogo Analiza, Novosibirsk, 29 (1976), pp.~3--10.

\bibitem{voigt1993}
{\sc M.~Voigt}, {\em List colourings of planar graphs}, Discrete Mathematics,
  120 (1993), pp.~215--219.

\end{thebibliography}

\end{document}